\documentclass[11pt,letterpaper]{amsart}
\usepackage{amssymb, amscd, amsfonts, euler, epic, xy, charter}
\usepackage{srcltx}

\def\QQ{{\mathbb Q}}

\def\p{\partial}

\def\half{\tfrac{1}{2}}
\def\G{{\Gamma}}

\def\arc{{\rm arc}}

\def\sep{{\rm sep}}

\def\bs{\backslash}

\def\Ccal{{\mathcal C}}

\def\Mcal{{\mathcal M}}

\def\Tcal{{\mathcal T}}

\newcommand\Lk{\operatorname{Link}}

\newcommand\Star{\operatorname{Star}}

\newtheorem{theorem}{Theorem}[section]
\newtheorem{lemma}[theorem]{Lemma}

\newtheorem{corollary}[theorem]{Corollary}

\theoremstyle{definition}
\newtheorem{definition}[theorem]{Definition}

\theoremstyle{remark} 
\newtheorem{remark}[theorem]{Remark}

\title{Connectivity of complexes of separating curves}
\author{Eduard Looijenga}
\email{E.J.N.Looijenga@uu.nl}
\address{Mathematisch Instituut, Universiteit Utrecht, P.O.~Box 80.010, NL-3508 TA Utrecht (Nederland)}
\thanks{Support by the Mathematical  Research Center at Tsinghua University is gratefully acknowledged.}



\begin{document}\sloppy
\maketitle
\hfill{\emph{In memory of Fritz Gr\"unewald (1949-2010)}

\begin{abstract}
We prove that the separating curve complex of a closed orientable surface of genus $g$ is
$(g-3)$-connected. We also obtain a connectivity property for a  separating curve complex 
of the open surface that is obtained by removing a finite set from a closed one, where 
it is assumed that the removed set is endowed with a partition and that the separating curves 
respect that partition. These connectivity statements   have implications for the algebraic topology of the moduli space of curves.
\end{abstract}

\section{Statements of the results}

Let $S$ be a connected oriented surface of genus $g$ with finite first Betti number $2g+n$ (i.e., a closed surface with $n$ punctures) and  make the customary assumption that $S$ has negative Euler characteristic: if $g=0$, then $n\ge 3$ and if $g=1$, then $n\ge 1$.
We recall that the \emph{curve complex} $\Ccal (S)$ of $S$ is the simplicial complex whose vertex set  consists of the isotopy classes of embedded (unoriented) circles in $S$ which do not bound in $S$ a disk or a cylinder. A finite set of vertices spans a simplex precisely when its elements can be represented by embedded circles that are pairwise disjoint. Thus, a closed 1-dimensional submanifold $A$ of $S$ with $k+1$ connected components such that every connected component of its complement has negative Euler characteristic defines a $k$-simplex $\sigma_A$ of $\Ccal (S)$ and every simplex of $\Ccal (S)$ is thus obtained.

This complex has proven to be quite useful in the study of the mapping class group of $S$. For the purposes of studying the Torelli group of $S$ a  subcomplex $\Ccal_\sep (S)$ of $\Ccal (S)$ can render a similar service. It is defined as the full subcomplex of  $\Ccal (S)$ spanned by the 
separating vertices of $\Ccal (S)$, where a vertex  is called \emph{separating} if a representative embedded circle separates $S$ into two components. Our main result for the case when  $S$ is closed is contained in the following theorem.

\begin{theorem}[$A_{g}$]\label{thm:Ag}
If $n\le 1$, then the simplicial complex $\Ccal_\sep (S)$ is  $(g-3)$-connected.
\end{theorem}

Previous work on this topic that we are aware  of concerns the case $n=0$. Farb and Ivanov announced
in 2005 \cite[Thm.\ 4]{fi} that $\Ccal_\sep (S)$ is connected for $g\ge 3$. Putman gave in \cite[Thm.\ 1.4]{putman:cn} another proof of this and showed that $\Ccal_\sep (S)$ is simply connected for $g\ge 4$ (\emph{op.\ cit.}, Thm.\ 1.11). In that paper he also mentions that Hatcher and Vogtmann have proved that $\Ccal_\sep (S)$ is $\lfloor \half(g-3)\rfloor$-connected for all $g$ (unpublished). 

\begin{remark}
Possibly the connectivity bound in Theorem \ref{thm:Ag}  is the best possible  for every positive genus. In a paper with Van der Kallen \cite{kl} we  showed that  the quotient of  $\Ccal_\sep (S)$ by the action of the Torelli group of $S$ has the homotopy type of a bouquet of $(g-2)$-spheres. In the situation of the theorem, $\Ccal_\sep (S)$ has dimension $2g-4+n$  so that its connectivity is half the dimension (as $n\le 1$). In particular, we cannot conclude that $\Ccal_\sep(S)$ is spherical.
\end{remark}

Before we state a version for the case $n\ge 2$, we point out a consequence that pertains to the moduli space of curves. Consider the Teichm\" uller space $\Tcal (S)$  of $S$, a contractible manifold  on which acts the mapping class group $\G (S)$. The action is proper and the orbit space may be identified with the moduli space $\Mcal_g$ of curves of genus $g$. The \emph{Harvey bordification} of $\Tcal (S)$, here denoted by  $\Tcal (S)^+\supset \Tcal (S)$, is a (noncompact) manifold with boundary with corners to which the action of $\G (S)$ naturally extends. This action is also proper and according to  \cite{looij} the orbit space $\Mcal_g^+:= \G(S)\bs \Tcal (S)^+$ is a compactification of $\Mcal_g$ that can also be obtained from the Deligne-Mumford compactification $\overline{\Mcal}_g\supset\Mcal_g$ as a `real oriented blowup' of its boundary $\Delta_g:= \overline{\Mcal}_g-\Mcal_g$.  The walls of $\Tcal (S)^+$ define a closed covering of the boundary $\p \Tcal(S)^+$ and any nonempty corner closure is an intersection of walls. As is well-known, the curve complex  $\Ccal (S)$ can be identified with the nerve of this covering of  $\p \Tcal(S)^+$. Since the corner closures  are contractible, Weil's nerve theorem implies that $\p \Tcal(S)^+$ has the same homotopy type as  $\Ccal (S)$. 

Let $\Delta_{g,0}\subset\Delta_g$ denote the irreducible component of the Deligne-Mumford boundary whose generic point parameterizes irreducible curves with one singular point. We may understand $\Mcal_g^c:=\overline{\Mcal}_g-\Delta_{g,0}$ as the moduli space of stable genus $g$ curves all of whose nodes are separating (which is equivalent to the irreducible components of the curve being smooth and with their genera summing up to $g$) and $\Delta^c_g:=\Delta_g-\Delta_{g,0}$ as the locus in  $\Mcal_g^c$ that parameterizes the singular ones among them. If  $\G$  is a subgroup of $\G (S)$ with the property that every Dehn twist along a separating curve in $S$  has a positive power lying in $\G$, then this defines a  (not necessarily finite) cover  $\tilde\Mcal^c_g\to \Mcal^c_g$.

\begin{corollary}
Suppose $\G\subset\G (S)$ is as above and  is in addition torsion free. If we denote by  $\tilde\Delta^c_g\subset \tilde\Mcal^c_g$ the preimage of $\Delta^c_g$, then the pair $(\tilde\Mcal^c_g,\tilde\Delta^c_g)$ is $(g-2)$-connected. Moreover, 
$H_k(\Mcal_g^c,\Delta^c_g;\QQ)=0$ for $k\le g-2$.
\end{corollary}
\begin{proof}
Let $\Tcal(S)^+_\sep$ be obtained from $\Tcal(S)^+$ by removing the walls that correspond to the nonseparating vertices of $\Ccal (S)$. Then $\Tcal(S)^+_\sep$ is the preimage of $\Mcal^c_g$ 
in $\Tcal(S)^+$. The same reasoning as above shows that $\p \Tcal(S)_\sep^+$ is homotopy equivalent to $\Ccal (S)_\sep$ and so $\p \Tcal(S)_\sep^+$ is $(g-3)$-connected. 
It follows that  we can construct a relative CW complex $(Z,\p \Tcal(S)_\sep^+)$ obtained from 
$\p \Tcal(S)_\sep^+$ by attaching cells of dimension $\ge g-1$ in a $\G(S)$-equivariant manner as to ensure that $Z$ is contractible and no nontrivial element of $\G (S)$  fixes a cell. Then $\G$ acts freely on $Z$ (as it does on the contractible space $\Tcal(S)^+_\sep$) and so there is a $\G$-equivariant homotopy equivalence  $Z\to \Tcal(S)^+_\sep$ relative to  $\p\Tcal(S)^+_\sep$.  It follows that  we also have a homotopy equivalence $\G\bs Z\to \tilde\Mcal^c_g$ relative to $\tilde \Delta^c_g$
and we conclude that   $(\tilde\Mcal^c_g,\tilde\Delta^c_g)$ is $(g-2)$-connected. 

The last assertion follows from the existence of a normal subgroup $\G\subset \G(S)$ of finite index that is torsion free. For if $\G$ is such a group, then 
$H_k(\Mcal_g^c,\Delta^c_g;\QQ)\cong H_k(\tilde\Mcal^c_g,\tilde\Delta^c_g;\QQ)^{\G(S)/\G}=0$ for $k\le g-2$.
\end{proof}

A similar statement holds for the universal curve $\Mcal_{g,1}$.

When $n>1$, we need to come to terms with the fact that the separability notion has no good heriditary properties: if $T$ is a closed surface, $A\subset T$ a compact 1-dimensional submanifold representing a simplex of $\Ccal (T)$ and $S$ a connected component  of $T-A$, then  a vertex of  $\Ccal (S)$ may split $S$, but not $T$. This happens precisely when the vertex in question separates two boundary components of $\p S$ that lie on the same  connected component of $T-S$. So the basic object should be, what Andy Putman  calls in  \cite{putman:cp}, a partitioned surface: a closed surface  minus a finite set, for which the removed set comes with a partition. This leads to the following definition.

\begin{definition}
Let $N$ be the set of points of $S$ at infinity (the cusps) and let  $P$ be  a partition of $N$. We call a vertex of $\Ccal (S)$ \emph{separating relative to $P$} if a representative embedded circle $\alpha\subset S$ has the property that $S-\alpha$ has two connected components each of which meets $N$ in a union of parts of $P$. We denote by $\Ccal (S,P)$ the full subcomplex of  $\Ccal (S)$ spanned by such vertices. 
\end{definition}

One might also understand $\Ccal (S,P)$ as the full subcomplex of  $\Ccal (S)$ spanned by the isotopy classes of embedded cycles  which are 
separating on  the surface $S^P$ that is obtained  by capping  off for  each part of $P$ the corresponding set of cusps by a sphere with that many holes. Notice that  $\Ccal (S,P)\subset\Ccal_\sep (S)$ and that we have equality when $P$ is discrete or $N$ is empty.

We shall prove Theorem \ref{thm:Ag} by induction and simultaneously with  

\begin{theorem}[$A_{g,n}$]\label{thm:Agn}
Suppose $g>0$ and $n=|N|>1$. Let $P$ be a partition of $N$. Then  $\Ccal (S,P)$ is $(g-2)$-connected.
\end{theorem}

To be precise, the induction starts with $g=0$, where the statements ($A_g$) and  ($A_{0,n}$) are trivially true and the induction strategy will be to show that 
\begin{itemize}
\item[(i)] ($A_{h,n}$) for $h<g$ implies ($A_g$) and
\item[(ii)] ($A_g$) and ($A_{h,k}$) for  $(h,k)<(g,n)$ (for the lexicographic ordering) imply ($A_{g,n}$).
\end{itemize}

I am indebted to Allen Hatcher for pointing out that the stronger version of Theorem \ref{thm:Agn} that I stated in a previous version was incorrect. Yet  it may be that some such statement might hold. For instance, if $r(P)$ denotes the number of nonempty parts of $P$ and $s(P)$ the number of parts with at least two elements, is it true that $\Ccal (S,P)$ is $(g+r(P)+s(P)-4)$-connected when $g>0$ (as I claimed in the  earlier version)? In case $g=0$, $\Ccal (S,P)$ is a complex of dimension $r(P)+s(P)-4$. Is this $(r(P)+s(P)-5)$-connected? In other words, is this complex spherical?
\\

I am grateful to the referee, whose meticulous job helped to improve the paper. The proof of Lemma \ref{lemma:htp}
follows a suggestion by the referee and simplifies my original one.

I also gratefully acknowledge support by the Mathematical Sciences Center of Tsinghua University at Beijing, where some of this work was done.

\section{Proofs}
Before we start off,  we  mention the following elementary fact that we will frequently use.

\begin{lemma}\label{lemma:join}
Let $X_i$ be a $d_i$-connected space ($d_i=-1$  means  $X_i\not=\emptyset$), where $i=1,\dots,k$. Then
the iterated join $X_1*\cdots * X_k$ is $(-2+\sum_{i=1}^k (d_i+2))$-connected.
\end{lemma}

\begin{proof}[Proof that ($A_{h,n}$) for $h<g$, all $n$, implies ($A_g$).]
So here $n\le 1$. We must show that $\Ccal_\sep (S)$ is $(g-3)$-connected. For $g<2$, there is nothing to show and so we may assume that $g\ge 2$. 
A theorem of Harer \cite[Thm.\ 1.2]{harer} asserts that $\Ccal (S)$ is  $(2g-3)$-connected. So it is certainly $(g-3)$-connected. 
Let  $\Ccal_k$ be the subcomplex  of $\Ccal (S)$ 
that is the union of $\Ccal_\sep (S)$ and the $k$-skeleton of $\Ccal(S)$. So $\Ccal_{-1}=\Ccal_\sep (S)$ and $\Ccal_k=\Ccal (S)$ for $k$ large.  
Notice that a finite set of vertices of $\Ccal (S)$ spans a simplex of $\Ccal_{k}$ if and only if no more than $k+1$ of these are nonseparating. Hence a minimal simplex of $\Ccal_k-\Ccal_{k-1}$ is represented by a compact 1-dimensional submanifold $A\subset S$  with $k+1$ connected components, each of which is nonseparating. We prove that the boundary of the star of such a simplex in $\Ccal_k$ is a  $(g-3)$-connected subcomplex of $\Ccal_{k-1}$. This property implies
that  the pair $(\Ccal(S),\Ccal_\sep(S))$ is $(g-2)$-connected. Since $\Ccal(S)$ is $(g-3)$-connected, it then follows that 
$\Ccal_\sep(S)$ is. 
Let $\{ S_i\}_{i\in I}$ be the set of connected components of $S-A$. Notice that if $g_i$ is the genus of $S_i$, then $g_i<g$. An Euler characteristic argument shows that  
\[
g-1= k+1+\sum_{i\in I} (g_i-1).
\]
We denote by $N_i$ the set of `cusps' of $S_i$, i.e.,  the finitely many points needed to make $S_i$ a closed surface. So an element  of $N_i$ is given by possibly a  cusp of $S$ (if it exists and if it is also a cusp of $S_i$) or by a connected component of $A$ in the  boundary of $S_i$ \emph{endowed with the orientation it receives as such}. The set $N_i$ comes with an evident partition $P_i$:  if $S$ has a cusp and $N_i$ contains it, then this cusp makes up a singleton part of  $P_i$ and 
any other two elements of $N_i$ belong to the same part of $P_i$ if and only if they come from connected components of $A$ that lie on the same connected component of $S-S_i$. (NB: beware that a connected component of $S-S_i$ could be simply a connected component $A_o$ of $A$; then its two orientations define a 2-element part of $P_i$.)
connected components of $A$ that bound $S_i$. Note that since the connected components of $A$ are nonseparating, we always have  $|N_i|\ge 2$. By our induction hypothesis 
$\Ccal (S_i,P_i)$ is then $(g_i-2)$-connected.
The boundary of the star of the $k$-simplex $\sigma_A$ defined by $A$ in $\Ccal_k$ lies in $\Ccal_{k-1}$ and can be identified with the $(|I|+1)$-fold join
\[
\p\sigma_A*  \big({*}_{\substack{i\in I}}\Ccal (S_i, P_i)\big).
\]
Since $\p\sigma_A$ is a combinatorial $(k-1)$-sphere, this join has  by Lemma \ref{lemma:join} connectivity at least 
$-2+k+\sum_{i\in I} g_i $. By the displayed formula above this is equal to $g-4+|I|$ and is therefore $\ge g-3$.
\end{proof}

The proof of $A_{g,n}$ begins with a discussion. We now assume that $g>0$ and $n\ge 2$.
Let $x\in N$. Notice that $S':= S\cup \{ x\}$ has  still negative Euler characteristic. We put $N':= N-\{ x\}$ and $P':=P|N'$. The goal is to compare  
$\Ccal (S',P')$ with $\Ccal (S,P)$. There is in general no forgetful map $\Ccal (S,P)\to \Ccal (S',P')$ because there will be vertices of $\Ccal (S,P)$ that do not give vertices of $\Ccal (S',P')$. Let us first identify this set of
vertices.

Denote by $\Sigma_x\subset N-\{ x\}$ the set of $y\in N-\{ x\}$  for which $\{ x,y\}$ is a union of parts of $P$. In other words, if 
$P_x$ denotes the part of $P$ that contains $x$, then $\Sigma_x$ is empty if $P_x$ has more than 2 elements, equals $P_x-\{ x\}$ if  $P_x$ is a 2-element set, and equals the set of $y\not=x$ for which $P_y$ is a singleton  in case $P_x=\{ x\}$. Then the  vertices of $\Ccal (S,P)$ that
have no image  in $\Ccal (S',P')$ are precisely the vertices $\alpha$ of $\Ccal_\sep (S)$ which for some $y\in\Sigma_x$ bound a disk neighborhood of $\{ x,y\}$ in $S\cup\{ x,y\}$.  Such a disk neighborhood can be thought of as a regular neighborhood of an arc in $S\cup\{ x,y\}$ connecting the two added cusps; this may help to explain why we have chosen to denote this set of vertices by  $\arc_{(S,P)}(x)$. Denote by $\Ccal (S,P)_x$ the full subcomplex of $\Ccal (S,P)$  spanned by the vertices not in $\arc_{(S,P)}(x)$. 

Observe that $\arc_{(S,P)}(x)$ is empty (so that $\Ccal (S,P)_x=\Ccal (S,P)$) if $\Sigma_x$ is. 
%(So  we can always find an $x\in N$ with nonempty $\Sigma_x$ unless $N$ is a singleton.)

\begin{lemma}\label{lemma:link}
The link in $\Ccal (S,P)$ of every vertex of $\arc_{(S,P)}(x)$ is a subcomplex of
$\Ccal (S,P)_x$ that projects isomorphically onto $\Ccal (S',P')$. 
\end{lemma}
\begin{proof}
A vertex of $\arc_{(S,P)}(x)$ defines a $y\in \Sigma_x$ and (up to isotopy) a closed disk $D$ in $S\cup\{x,y\}$ that is a neighborhood of  $\{ x,y\}$.  The inclusion $S-D\subset S'$  identifies the link in question with $\Ccal (S',P')$.
\end{proof}

Denote by $\tilde P$ the  refinement of $P$ which coincides with $P$ on $N-P_x$ and 
partitions $P_x$ further into $\{ x\}$ and $P_x-\{ x\}$. So $\tilde P'=P'$.
It is clear that $\Ccal (S,P)$ is a subcomplex of $\Ccal (S,\tilde P)$. Notice that $\arc_{(S,P)}(x)=\Ccal (S,P)\cap\arc_{(S,\tilde P)}(x)$ (we have $\arc_{(S,P)}(x)=\arc_{(S,\tilde P)}(x)$ unless $|P_x|=2$) and 
$\Ccal (S,P)_x=\Ccal (S,P)\cap\Ccal (S,\tilde P)_x$. We denote by $f$  the forgetful simplicial map $\Ccal (S,\tilde P)_x\to \Ccal (S',P')$ so that we have the diagram
\[
%\begin{CD}\[
\begin{array}{ccccc}
C(S,P)_x  & \subset  &  C(S,\tilde P)_x  & {\buildrel f\over\longrightarrow}&  C(S',P'). \\
\cap  &   & \cap &  &\\
C(S,P)  & \subset   & C(S,\tilde P)  & & 
\end{array}
\]

\begin{lemma}\label{lemma:htp}
The map $f: \Ccal (S,\tilde P)_x\to \Ccal (S',P')$ is a homotopy equivalence.
\end{lemma}
\begin{proof}
Choose  an arc $\gamma$ which connects $x$ with another point of $N$ and defines a vertex of $\arc_{(S,P)}(x)$ and 
observe that the full subcomplex  $K\subset\Ccal (S, \tilde P)_x$ spanned by vertices that avoid $\gamma$ defines a section of $f$ (the inclusion $S-\gamma\subset S'$ is isotopic to a homeomorphism). We shall prove that  $\Ccal (S, \tilde P)_x$ admits $K$ as a deformation retract. (The proof will in fact show that each fiber of $|f|$ is a tree and essentially produces for every element of $|\Ccal (S, \tilde P)_x|$ the unique path in its $|f|$-fibre that connects it to the point of $|K|$.) 

Denote by $K_r\subset\Ccal (S, \tilde P)_x$ the subcomplex whose simplices can be represented by a closed submanifold $A\subset S$ which meets $\gamma$ transversally in at most $r$ points. This defines a filtration $K=K_0\subset K_1\subset K_2\subset\cdots$ whose union is
$\Ccal (S, \tilde P)_x$. Although this filtration is infinite, it is enough to construct for every $r\ge 0$ a  deformation retraction of $|K_{r+1}|$ onto $|K_r|$, for in the simplicial setting an infinite sequence of deformation retractions still gives a deformation retraction.

We do this per simplex:
if  $\sigma$ is a simplex of $K_{r+1}$ that is not in $K_r$ and is minimal for this property, then its link in $K_{r+1}$ lies in $K_r$ and so it suffices to define for such a $\sigma$ a  deformation retraction $h_\sigma$ of $|\Star_{K_{r+1}}(\sigma)|$ onto $|\Lk_{K_{r+1}}(\sigma)|$.

The simplex $\sigma$  is represented  by a closed submanifold  $A\subset S$ of which every connected component meets $\gamma$ transversally and is  such that
$A\cap\gamma$ has cardinality $r+1$ (a number that cannot be made smaller in its isotopy class). 
Let $x_0$ be the point of $A\cap\gamma$ closest to $x$. Denote by $\alpha_0$ the connected component of $A$ which contains $x_0$ and choose in $S'$ a thin regular neighborhood of the union of $\alpha_0$ and the 
subarc  of $\gamma$ which connects $x_0$ with $x$. The boundary of that neighborhood  has two connected components. Both lie in $S$ and only one of them is isotopic to $\alpha_0$. Denote by  $\alpha'_0$ the other boundary component. If $\tau$ is a simplex of $K_{r+1}$ which contains $\sigma$, then adding $\alpha'_0$ to $\tau$ gives also a simplex $\tau'$ of 
$K_{r+1}$ and the  codimension one face  $\tau''$ of $\tau'$ obtained by removing $\alpha_0$ is contained in $K_r$.  So if we  regard $|\Star_{K_{r+1}}(\sigma)|$ as the cone over $|\Lk_{K_{r+1}}(\sigma)|$ with  the barycenter of $\sigma$ as its vertex, then there is a simplicial map from this cone to its base which  sends the barycenter to $\alpha'_0$ and is the identity on the base. Its  geometric realization yields the desired $h_\sigma$.
\end{proof}

\begin{corollary}\label{cor:connect}
The complex $\Ccal (S,P)\cup \Ccal (S,\tilde P)_x$ is canonically homotopy equivalent to the join 
$\arc_{(S,P)}(x)\ast\Ccal (S',P')$ (where $\arc_{(S,P)}(x)$ is  discrete).
\end{corollary}
\begin{proof}
The set of vertices of $\Ccal (S,P)\cup \Ccal (S,\tilde P)_x$ not in $\Ccal (S,\tilde P)_x$ is  $\arc_{(S,P)}(x)$. The link of any such vertex is contained in $\Ccal (S,\tilde P)_x$ and by  Lemma \ref{lemma:link} that link projects isomorphically onto $\Ccal (S',P')$. 
In view of Lemma \ref{lemma:htp} this implies that the inclusion of this link in $\Ccal (S,\tilde P)_x$ is also a homotopy equivalence. Hence the natural inclusion $\Ccal (S,P)\cup \Ccal (S,\tilde P)_x\subset \arc_{(S,P)}(x)\ast  \Ccal (S,\tilde P)_x$ is a homotopy equivalence. The corollary follows.
\end{proof}

\emph{From now on we assume that $A_g$ holds and that $A_{h,k}$ holds  for all $(h,k)$ smaller than $(g,n)$ for the lexicographic ordering. Our goal is to prove $A_{g,n}$.} 

\begin{lemma}\label{lemma:connect}
The pair 
$(\Ccal (S,P)\cup \Ccal (S,\tilde P)_x,\Ccal (S,P))$ is $(g-1)$-connected.
\end{lemma}
\begin{proof} 
If $P_x=\{ x\}$, then $\tilde P= P$ and there is nothing to show.
We therefore assume that $P_x$ has more than one element.
Denote by $\Ccal_k$  the subcomplex 
of $\Ccal (S,P)\cup\Ccal (S,\tilde P)_x$ that is the union of $\Ccal (S,P)$ and the $k$-skeleton of $\Ccal (S,P)\cup\Ccal (S,\tilde P)_x$: a finite set of vertices of $\Ccal (S,P)\cup\Ccal (S,\tilde P)_x$ spans a simplex of $\Ccal_{k}$ if and only if no more than $k+1$ of these separate $x$ from $P_x-\{x\}$. Notice  that $\Ccal_{-1}=\Ccal (S,P)$ and  $\Ccal_k=\Ccal (S,P)\cup\Ccal (S,\tilde P)_x$ for $k$ large. A minimal simplex of $\Ccal_k-\Ccal_{k-1}$ is represented by a compact 1-dimensional submanifold $A\subset S$  with $k+1$ connected components, each of which separates $x$ from $P_x-\{ x\}$ (the graph that is associated to $A$ is then a string with $k+2$ nodes). We prove that the boundary of the star of such a simplex in $\Ccal_{k-1}$ is
$(g-2)$-connected if $g>0$. 
We enumerate the connected components of $A$ as $\alpha_0,\dots ,\alpha_k$ and the connected components of $S-A$ as $S_0,\dots ,S_{k+1}$ such that $\alpha_i$ is a boundary component of $S_i$ and $S_{i+1}$ and so that  $S_0$ resp.\ $S_{k+1}$ is  punctured by $x$ resp.\  $P_x-\{ x\}$. The cusps of $S-A$  are naturally indexed by $\hat N:= N\sqcup\{ i_\pm\}_{i=0}^k$, where $i_-$ resp.\ $i_+$ corresponds to the cusp defined by $\alpha_i$ on $S_i$ resp.\ $S_{i+1}$. Let $\hat N_i\subset \hat N$ index the set of cusps on $S_i$. Denote by $P_i$ the partition of $(N-P_x)\cap S_i$ that is simply the restriction of $P$ and denote by $\hat P_i$ the partition of $\hat N_i$ that on $(N-P_x)\cap S_i$ is equal to $P_i$ and has what remains of $\hat N_i$ as a single part. So this new part is 
$ \{ x ,0_+\} $ for $i=0$,  $\{ (i-1)_-,i_+\}$ for $0<i<k+1$ and $(P_x-\{x\})\cup\{ k_-\}$ for $i=k+1$.

The reason for introducing these partitions  is that we can now observe that the boundary of the star of the $k$-simplex $\sigma_A$ defined by $A$ in $\Ccal_k$ lies in $\Ccal_{k-1}$ and can be identified with the iterated join
\[
\p \sigma_A * \Ccal (S_0,\hat P_0)*\cdots * \Ccal (S_{k+1},\hat P_{k+1}).
\]
It is then enough to show that this join is $(g-2)$-connected for $g>0$. 
Since $\p \sigma_A$ is a $(k-1)$-sphere, it is $(k-2)$-connected.
The connectivity of a factor $\Ccal (S_i,\hat P_i)$ with $g_i>0$ is at least $g_i-2$. So by Lemma \ref{lemma:join} the connectivity of the above join is at least $-2+k+\sum_{\{i: g_i>0\}} g_i= g+k-2\ge g-2$.
\end{proof}

\begin{proof}[Proof of ($A_{g,n}$)]
We must show that $\Ccal (S,P)$ is $(g-2)$-connected.
In view of Lemma \ref{lemma:connect} it suffices to show that $\Ccal (S,P)\cup \Ccal (S,\tilde P)_x$ has that property.

If $\arc_{(S,P)}(x)=\emptyset$, then $n>2$ and so our induction hypothesis implies that $\Ccal (S',P')$ is  $(g-2)$-connected by $A_{g,n-1}$. It follows from Corollary \ref{cor:connect} that $\Ccal (S,P)\cup \Ccal (S,\tilde P)_x$ is homotopy equivalent to
$\Ccal (S',P')$ and hence is $(g-2)$-connected. 

If $\arc_{(S,P)}(x)\not=\emptyset$, then we may have $n=2$. At least we know that $\Ccal (S',P')$ is $(g-3)$-connected
(invoke $A_g$ if $n=2$).
But since $\Ccal (S,P)\cup \Ccal (S,\tilde P)_x$ is homotopy equivalent to
$\arc_{(S,P)}(x)\ast \Ccal (S',P')$ (by Corollary \ref{cor:connect}), it is $(g-2)$-connected. 
\end{proof}

\end{document}